\documentclass[12pt]{article}
\usepackage{amsmath}
\usepackage{amsthm}
\usepackage{amsfonts}

\newtheorem{thm}{Theorem}[section]
\newtheorem{lemma}[thm]{Lemma}

\newtheorem{que}[thm]{Question}

\newtheorem{con}[thm]{Conjecture}
\newtheorem{rem}[thm]{Remark}
\newcommand{\R}{{\mathbb R}}
\newcommand{\Z}{{\mathbb Z}}

\def\Cox{\hfill{QED}}

\def\st{:\,}

\def\P{{\bf P}}

\begin{document}

\title{On the effect of adding $\epsilon$-Bernoulli percolation to everywhere
percolating subgraphs of $\Z^d$}
\author{Itai Benjamini, Olle H\"aggstr\"om and Oded Schramm}

\maketitle

\begin{abstract}
We show that adding $\epsilon$-Bernoulli percolation  
to an everywhere percolating 
subgraph of $\Z^2$ results in a graph which has large scale 
geometry similar to that of supercritical Bernoulli percolation, in
various specific senses. 
We conjecture similar behavior in higher dimensions.
\end{abstract}

\section{Introduction}

A subset $X$ of the edges of the standard $d$-dimensional cubic lattice 
$\Z^d$ is said to be {\bf percolating everywhere} if every vertex 
of $\Z^d$ is contained in an infinite connected component of $X$.
Examples of such subgraphs are foliations by lines, and spanning forests.
In this note we study the effect of adding small noise to 
the geometry of such
subgraphs of $\Z^d$. We will argue that if $X$
is percolating everywhere, then adding $\epsilon$-Bernoulli
percolation acts as a unifying 
operation on the geometric structure of the subgraph; see 
Conjecture \ref{mainc} and Theorem \ref{maint} below. 
By ``adding $\epsilon$-Bernoulli percolation'', we mean that each edge
that is not in $X$, is added independently with probability $\epsilon$. 

So far, we can prove our claims only in dimension two. Our proofs make
crucial use of planar duality, so that new ideas clearly are needed to
make progress in higher dimensions.

Here is some motivation for our study. 
By $p$-Bernoulli percolation on an infinite graph $G$,
we mean the usual bond percolation process, where each edge is removed with
probability $1-p$ and kept with probability $p$. By $p_c(G)$, we denote the
infimum over all $p\in[0,1]$ such that $p$-Bernoulli percolation on $G$
has infinite clusters with positive probability.
 An outstanding open problem
in percolation theory (see e.g.\ Grimmett (1989)) is to determine whether
at criticality $p=p_c$ there are infinite clusters; the answer is believed
to be no for all $d\geq 2$. Meditating over this problem, 
one is naturally lead to a search for
conditions on $X\subset \Z^d$ which guarantee $p_c(X) < 1$.
If it could be shown that infinite Bernoulli-percolation clusters $W$ satisfy
$p_c(W)<1$, then it would follow that 
there are no infinite clusters at $p=p_c$.
In particular, a negative answer to the following question would
answer the problem of the existence of infinite clusters at $p_c$.

\begin{que}
\label{q1}
Is there an invariant finite energy percolation $X$ on $\Z^d$, which  
a.s.\ percolates and satisfies $p_c(X)=1$?
\end{que}

An invariant percolation is a random subgraph of $\Z^d$ whose distribution
is invariant under the automorphisms of $\Z^d$.
Finite energy percolation was first considered by Newman and
Schulman (1982), and is the same as 
deletion and insertion tolerance in the sense of
Lyons and Schramm (1999): deletion (resp.\ insertion) tolerance means that 
the conditional probability that an edge is absent (resp.\ present) given 
the status of all other edges is strictly positive. 
One way of constructing examples of insertion 
tolerance percolation is to
add independent $\epsilon$-Bernoulli percolation to any given percolation
process. In Section 3, we will give an 
an example of an invariant insertion tolerant percolation process $X$
obtained via adding  
$\epsilon$-Bernoulli percolation, which percolates but for which $p_c(X)=1$. 
In that example, large chunks of
vertices in $\Z^d$ are in finite connected components of the percolation.
This observation led us to 
% Conjecture \ref{mainc} below. 

% Our vague conjecture is that if $Y = Y(X, \epsilon)$ is obtained by adding 
% $\epsilon$-Bernoulli percolation to an everywhere percolating subgraph $X$,
% then a.s.\ $Y$ ``is very similar'' to supercritical Bernoulli percolation. 
% To make sense of this, we mean in particular
% the following list of properties.

\begin{con}
\label{mainc}
Let $X$ be a fixed
everywhere percolating subgraph of $\Z^d$, and let $Y=Y(X,\epsilon)$ be
obtained from $X$ by adding $\epsilon$-Bernoulli percolation.
For any $\epsilon > 0$, we have
\begin{description}
\item{\rm (i)} $Y$ is connected a.s.
\item{\rm (ii)} $p_c(Y) <1$ a.s.
\item{\rm (iii)} $Y$ percolates in the upper half-space a.s.
\item{\rm (iv)} 
A renormalized version of $Y$ dominates supercritical Bernoulli 
percolation.
\end{description}
\end{con}

\begin{thm}
\label{maint}
In dimension $d=2$, with $X$, $\epsilon$ and $Y$ as above, properties
{\rm (i)}, {\rm (ii)}, {\rm (iii)} and {\rm (iv)} hold.  
\end{thm}

We need to explain what is meant by the renormalization in item (iv). For
a positive integer $n$ and a vertex $x\in \Z^d$, let $\Lambda(x,n)$ denote
the box $x+ [-\frac{n}{2}, \frac{n}{2}]^d$ of side-length $n$ centered at
$x$. If $x$ and $y$ are nearest neighbors in $\Z^d$, then the vertices
$nx$ and $ny$ are said to be {\bf closely connected} (in $Y$) if there is
a path in $Y$ from $nx$ to $ny$ inside $\Lambda(nx,n) \cup \Lambda (ny, n)$.
A renormalized version $\tilde{Y}_n$ of $Y$ is defined as the percolation
in $\Z^d$ where each edge $\langle x,y \rangle$ is included in
$\tilde{Y}_n$ if and only if $nx$ and $ny$ are closely connected in $Y$. 
Property (iv) then says that there
exist $p>p_c(\Z^d)$ and $n$, such that $\tilde{Y}_n$ stochastically 
dominates $p$-Bernoulli percolation on $\Z^d$. 

Our proof of Theorem \ref{maint} (iv) will in fact show the stronger result
that for any $p<1$, $\tilde{Y}_n$ dominates $p$-Bernoulli percolation for
all sufficiently large $n$.

\begin{rem}
\label{rem1}
{\rm 
If $X$ is an everywhere percolating realization of some 
invariant percolation on $\Z^d$, then  a.s.\ property (i) holds
for $Y=Y(X,\epsilon)$,
by the encounter points argument of Burton and Keane (1989). }
\end{rem}

\begin{rem}
{\rm 
Say that a subgraph
$X$ of $\Z^d$ is densely percolating, if there is some $R >0$ such that
any ball of radius $R$ in $\Z^d$ intersects an infinite connected component of
$X$. A straightforward extension of our arguments show that an analogue of 
Theorem \ref{maint} holds for densely percolating subsets of $\Z^2$ (note 
that property (i) of course has to be replaced by uniqueness of the infinite
cluster, and the definition of renormalization in (iv) has to be
modified slightly to allow e.g.\ the point $nx$ to be replaced by some
percolating point in its $R$-neighborhood). }
\end{rem}

\section{Proofs}

A main ingredient in our proofs is the use of planar duality. For a
(possibly random) edge configuration $X$ in $\Z^2$, let $X^*$ denote
the edge configuration in the planar dual $\Z^2_{dual}$
of $\Z^2$, where each edge in $\Z^2_{dual}$
is present if and only if the (unique) edge in $\Z^2$ that
crosses it is absent from $X$. 

\medskip\noindent
{\bf Proof of Theorem \ref{maint} (i):} If $X$ is percolating everywhere,
then it contains no finite connected components, so that the dual $X^*$ 
contains no circuits. Hence, 
for any fixed $x,y \in \Z^2_{dual}$, there is at most
one self-avoiding path in $X^*$ connecting them. This path has, of course,
length at least $|x-y|_1$, where $|\cdot|_1$ denotes $L^1$-distance in
$\R^2$. 

That $Y$ is obtained from $X$ via $\epsilon$-Bernoulli addition of edges, is
the same as saying that $Y^*$ is obtained from $X^*$ by randomly deleting each
edge in $X^*$ independently with probability $\epsilon$. Letting 
$\stackrel{Y^*}{\longleftrightarrow}$ denote connectivity in the $Y^*$
configuration, we get for any $x,y\in \Z^2_{dual}$ that
\begin{equation} \label{eq:pairwise_connected_dual}
\P(x\stackrel{Y^*}{\longleftrightarrow}y) \leq (1-\epsilon)^{|x-y|_1} \, .
\end{equation}
For any $x\in \Z^2_{dual}$ and any $k \geq 1$, there are exactly $4k$
vertices in $\Z^2_{dual}$ at $L^1$-distance $k$ from $x$. Summing 
(\ref{eq:pairwise_connected_dual}) over all $y\in \Z^2$, we get that
the expected number of vertices that are connected to $x$ in $Y^*$ is at most
\begin{equation} \label{eq:pairwise_connected_dual_sum}
1+ 4\sum_{k=1}^\infty k(1-\epsilon)^k \, < \, \infty \, .
\end{equation}
Hence the connected component of $Y^*$ containing $x$ is finite a.s., and
$Y^*$ is therefore a.s.\ a forest of finite trees. This implies that $Y$ is
connected a.s.
\Cox

\medskip\noindent
Our next task will be to prove Theorem \ref{maint} (iv); once this is done, 
properties (ii) and (iii) will be simple corollaries. For the proof
of (iv), the following lemma is useful.

\begin{lemma} \label{lem:path_radius}
For any nearest neighbors $x$ and $y$ in $\Z^2$, let $E_k^{x,y}$ 
denote the event that $x$ and $y$ are {\bf not} connected by any
path in $Y$ that is contained in the box $\Lambda(x,k)$. There exists a
constant $c>0$ (depending only on $\epsilon$) such that
\[
\P(E_k^{x,y}) \leq e^{-ck}
\]
for all $k$.
\end{lemma} 

\noindent
{\bf Proof:}
Here is a particular way of finding a path in $Y$ from $x$ to $y$: If the
edge $\langle x,y \rangle$ is present in $Y$, then use that edge. If
that edge is not present, then the corresponding edge
$\langle x,y \rangle^*$ is present in $Y^*$. We can then find a path from
$x$ to $y$ in $Y$ by going around the $Y^*$-component $T^*_{x,y}$ containing 
$\langle x,y \rangle^*$ clockwise, following the outer boundary
of $T^*(x,y)$. 
If $T^*_{x,y}$ is contained in $\Lambda(x, k-1)$, then the path
we just constructed is
contained in $\Lambda(x,k)$. By inspecting the summands in
(\ref{eq:pairwise_connected_dual_sum}), 
we see that the probability that $T^*_{x,y}$ is 
{\em not} contained in $\Lambda(x, k-1)$ decays exponentially in $k$,
which is what we needed. 
\Cox

\medskip\noindent
{\bf Proof of Theorem \ref{maint} (iv):}
Let $x$ and $y$ be nearest neighbors in $\Z^2$, and let $A^{x,y}_n$ 
be the event that $xn$ and $yn$ are closely connected. Let 
$z_0=x, z_1, z_2, \ldots, z_{n-1}, z_{n}=y$ be the vertices on the unique
shortest path from $nx$ to $ny$ in $\Z^2$. Clearly, 
\[
A_n \supset \neg\left( \cup_{i=0}^{n-1}E_n^{x_{i},x_{i+1}}\right)
\]
so that
\begin{eqnarray*}
\P(\neg A_n) & \leq & \P \left( \cup_{i=0}^{n-1}E_n^{x_{i},x_{i+1}}\right) \\
& \leq & \sum _{i=0}^{n-1} \P(E_n^{x_{i},x_{i+1}}) \\
& \leq & ne^{-cn}
\end{eqnarray*}
(where $c$ is as in Lemma \ref{lem:path_radius}). Hence
$\P(A_n) \geq 1-ne^{-cn}$, which tends to $1$ as $n\rightarrow \infty$.
Therefore, the probability that an edge in the renormalized process
$\tilde{Y}_n$ is present tends to $1$ as $n$ tends to infinity. This
observation
does not immediately imply the desired stochastic domination, because
the edges do not appear in $\tilde{Y}_n$ independently.

However, $\tilde{Y}_n$ is easily seen to be
a {\bf 1-dependent} percolation process, meaning the following: if 
$B_1, B_2 \subset \Z^2$ are two disjoint edge sets where no edge in
$B_1$ shares an endpoint with an edge in $B_2$, then $\tilde{Y}_n(B_1)$ and
$\tilde{Y}_n(B_2)$ are independent (this is simply because $\tilde{Y}_n(B_1)$ 
and $\tilde{Y}_n(B_2)$ depend on disjoint edge sets in $Y$). Theorem 6.5 of
Liggett, Schonmann and Stacey (1997) tells us that for any $p<1$, we can
find a $p'<1$ such that any 1-dependent percolation processes with
edge marginals greater than $p'$ dominates $p$-Bernoulli percolation.
So now we only need to pick $p\in(p_c(\Z^2),1)$, then pick $p'$ as in the
Liggett--Schonmann--Stacey theorem, and finally pick $n$ large enough to
guarantee that the edge marginals in $\tilde{Y}_n$ are greater than $p'$.
\Cox

\medskip\noindent
{\bf Proof of Theorem \ref{maint} (ii):}
Pick $n$ large enough so that property (iv) holds, i.e. so that
$\tilde{Y}_n$ dominates $p$-Bernoulli percolation for some $p>p_c(\Z^2)$. 
For $q\in(0,1)$, let 
$W_q$ be an independent $q$-Bernoulli percolation on $\Z^2$,
so that $Y\cap W_q$ is a $q$-Bernoulli percolation on $Y$. 
Pick $q$ close enough to $1$, so that for any $x$
the probability that there is an edge in $\Lambda(x,n)\setminus W_q$
is at most $(p-p_c(\Z^2))/4$. If we now 
thin $\tilde{Y}_n$ by removing any edge $\langle x,y\rangle$ in $\tilde{Y}_n$
such that some edge in $\Lambda(nx, n) \cup \Lambda(ny,n)$ 
is not in $W_q$,
 then the thinned $\tilde{Y}_n$-process
dominates Bernoulli percolation with parameter $(p+p_c(\Z^2))/2$, so
it still percolates. But if the thinned $\tilde{Y}_n$-process percolates,
then, clearly, so does $Y\cap W_q$.
\Cox

\medskip\noindent
{\bf Proof of Theorem \ref{maint} (iii):} This is immediate from
property (iv) and the fact that supercritical Bernoulli percolation
on $\Z^2$ percolates also in the upper half-plane; the latter result can be
found e.g.\ in Kesten (1982). 
\Cox

\section{An example}

We finally present an example of an invariant
percolation $X\subset\Z^d$ ($d\geq 2$), which has infinite clusters,
and nevertheless also has the property that
for any $\epsilon\in(0, p_c(\Z^d))$, adding $\epsilon$-Bernoulli
percolation to $X$ a.s.\ produces a graph $Y=Y(X,\epsilon)$ with
$p_c(Y)=1$. Note that $Y(X,\epsilon)$ is
insertion tolerant when $\epsilon>0$. 

Vaguely speaking, $X$ will be constructed by taking the full configuration
(all edges present), and removing edges from large annuli (of drastically
different sizes) in such a way that the outside and the inside connect only by
a thin thread. The annuli are spread out randomly, in such a way that the
origin is a.s.\ surrounded by infinitely many of them. $X$ then percolates, 
but the threads are cut when doing Bernoulli-thinning of $X$, and adding 
$\epsilon$-Bernoulli percolation doesn't help in bridging the annuli. 

The precise construction of $X$ is as follows.
Consider independent random
variables $\bigl\{a(x,n)\st (x,n)\in\Z^d\times\{1,2,\dots\}\bigr\}$
where 
\[
\P\bigl(a(x,n)=1\bigr)=2^{-dn}=1-\P\bigl(a(x,n)=0\bigr) \, .
\]
Let 
\[
W(x,n)=\Lambda(x,2^n)\setminus\Lambda(x,2^n-2^{n/2}) \, ,
\]
where, as before, $\Lambda(x,n)=x+ [-\frac{n}{2}, \frac{n}{2}]^d$.
Let $b(x,n)$ be the indicator of the event that
$a(y,k)=0$ for every $(y,k)\neq (x,n)$ such that
$n\geq k$ and $W(y,k)\cap W(x,n)\neq\emptyset$.
Let $W'(x,n)$ be the set of edges of the grid $\Z^d$ which
are inside $W(x,n)$, except those on the line
$x+\R\times\{0\}\times\cdots\times\{0\}$.
Finally, let $X$ consist of all edges of $\Z^2$ that are {\em not} in 
the set
\begin{equation}
\label{e.Xdef}
% \bigcup_{x\in \Z^d \atop n\in\{1,2,\ldots\}} \Bigl\{W'(x,n)\st
% (the range of union is inside)
\bigcup \Bigl\{W'(x,n)\st
(x,n)\in\Z^d\times\{1,2,\dots\},\ a(x,n)=b(x,n)=1\Bigr\}.
\end{equation}
It is immediate that a.s.\ $X$ has an infinite connected
component, and it is also straightforward to verify that
if in (\ref{e.Xdef}) we replace $W'$ with $W$,
then a.s.\ no infinite cluster remains.
Furthermore, using the (well-known, see e.g.\ Grimmett (1989)) exponential
tail of the cluster size distribution for subcritical Bernoulli percolation
on $\Z^d$, we see that for $\epsilon< p_c(\Z^d)$ the probability of 
bridging an
annulus $W(x,n)$ tends to $0$ as $n\rightarrow \infty$. It follows easily
that 
$p_c\bigl(Y(X, \epsilon)\bigr)=1$ for any $\epsilon< p_c(\Z^d)$.

% \noindent
% {\bf Aknowledgement:} Thanks to all the participants in the Yehuda Hills 
% picnic, which provided the right atmosphere leading to this work.

\bigskip

\noindent
\texttt{itai@wisdom.weizmann.ac.il} \\
\texttt{http://www.wisdom.weizmann.ac.il/\string~itai/}
\medskip

\noindent
\texttt{olleh@math.chalmers.se} \\
\texttt{http://www.math.chalmers.se/\string~olleh/}
\medskip

\noindent
\texttt{schramm@wisdom.weizmann.ac.il} \\
\texttt{http://www.wisdom.weizmann.ac.il/\string~schramm/}


\begin{thebibliography} {abc}

\bibitem{bk}
Burton, R.M. and Keane, M.S. (1989) Density and uniqueness in percolation,
{\it Commun. Math. Phys.} {\bf 121}, 501--505.

\bibitem{G}
Grimmett, G.R. (1989) {\em Percolation}, Springer, New York.

\bibitem{K} 
Kesten, H. (1982) {\em Percolation Theory for Mathematicians}, Birkh\"auser,
Boston. 

\bibitem{lss}
Liggett, T., Schonmann, R. and Stacey, A. (1997)
Domination by product measures. {\it Ann. Probab.} {\bf 25}, 71--95. 

\bibitem{ls}
Lyons, R. and Schramm, O. (1999) Indistinguishability of percolation clusters,
{\it Ann. Probab.}, to appear.

\bibitem{NS}
Newman, C.M. and Schulman, L.S. (1982) Infinite clusters in percolation
models, {\em J. Statist. Phys.} {\bf 26}, 613--628. 

\end{thebibliography}
\end{document}